\documentclass[11pt]{amsart}
\usepackage{amssymb}
\usepackage{amscd}
\usepackage{eucal}
\numberwithin{equation}{section}

\begin{document}

\dedicatory{An expanded version of a talk given at the AMS
Session\\ 
honoring George Mackey\\ New Orleans, January 7, 2007}

\title[Borel Classification Theory]{Classifying the Unclassifiables}
\author{Edward G. Effros}
\address{Department of Mathematics\\UCLA, Los Angeles, CA 90095-1555}
\email[Edward G. Effros]{ege@math.ucla.edu}
\date{March 11, 2007}
\maketitle
\section{Introduction}

George Mackey devoted his career to exploring the basic unity of
mathematics, and to understanding its relationships to modern physics. He was a major figure in the flowering of group
representation theory and related portions of functional analysis that occured in the middle decades of the last
century. Sadly, he, along with many of his colleagues who participated in this endeavor, is no longer with us.

Borel classifiability was one of Mackey's most innovative ideas. Although he only formulated it for unitary group representations, his theory can be adapted to a wide range of mathematical problems. Owing to the efforts of descriptive set-theorists, there has been remarkable progress in this field in the last twenty years. I have sketched some of of these developments, and I have suggested some open problems. I conclude with a few remarks about Mackey's attempts to come to terms with quantum mechanics. 

For a more complete and balanced view of George Mackey's career, the reader should turn to Varadarajan's masterful article.

A number of individuals have made very helpful suggestions for improving earlier versions of this paper. I am particularly indebted to Aleko Kechris, Richard Kadison, and Marc Rieffel.

\section{Type I if and only if smooth dual}
Functional analysis provides the essential tools for generalizing linear algebra to infinite dimensional spaces. Mackey's earliest significant work was in the area of locally convex spaces. Under the influence of such pioneers
as von Neumann, Stone, and Weyl, Mackey shifted his research to the study of infinite dimensional unitary
representations of locally compact groups. In order to generalize the finite dimensional techniques of classical representation theory, it was necessary for him
to replace \emph{direct sums} of finite dimensional spaces by von Neumann's theory of \emph{direct integrals} of Hilbert spaces \cite{vn}. In an even more striking
departure from the purely algebraic context, he used Murray and von
Neumann's notion of \emph{continuous dimensions} \cite{vn2} to compare invariant
subspaces. Thus Mackey's theory is a synthesis of sophisticated
measure-theoretic arguments and operator algebraic notions with the elegant machinery of the group theorists. Perhaps the capstone of this work was his far-reaching theory of induced representations (the ``Mackey machine''), which enables one to relate infinite dimensional group representations of a group to those of its subgroups.  

In simplistic terms, Mackey's purpose was to classify the unitary
representations of locally compact groups. As in the algebraic theory, the
first step is to decompose a representation into irreducibles. A finite
dimensional unitary representation $\pi $ of a finite group $G$ is unitarily
equivalent to a direct sum of irreducible representations 
\begin{eqnarray}
\pi \simeq \pi _{1}\oplus \ldots \oplus \pi _{r}.
\end{eqnarray}
One may collect equivalent irreducibles, and thus obtain a unique
decomposition 
\[
\pi \simeq n_{1}\theta _{1}\oplus \cdots \oplus n_{s}\theta _{s},
\]
where the $\theta _{j}$ are inequivalent irreducibles, and $n\theta =\theta
\oplus \cdots \oplus \theta $, with $n$ denoting the ``multiplicity'' of $\theta .
$ Letting $\hat{G}$ be the \emph{dual }of $G,$ i.e., the set of all unitary
equivalence classes of irrreducible unitary representations of $G,$ we may
rewrite this in the canonical form 
\begin{eqnarray}
\pi \simeq \sum\nolimits_{\theta \in \hat{G}}^{\oplus }n_{\theta }\theta
\end{eqnarray}
where $n_{\theta }\in \Bbb{N}\cup \left\{ 0\right\} .$ If all of the $%
n_{\theta }$ are equal to 0 or 1, then $\pi $ is said to be multiplicity
free, whereas if $\pi=n\theta$ for a single irreducible representation, $\pi$ is said to be a factor representation. Equivalently, the commutant operator algebra $\pi (G)^{\prime }$ is commutative (respectively, has scalar center, i.e., it is a factor).

Turning to the infinite dimensional unitary representations of a (second countable)
locally compact group $G$, Mautner (\cite{Mau}, Th. 1.2 -- see also \cite{Godem}, Th.7) generalized (1.1) by
showing that any representation $\pi $ has a direct integral decomposition
into irreducible representations 
\[
\pi \simeq \int_{X}^{\oplus }\pi _{x}d\mu (x),
\]
where $X$ is a standard Borel space. We recall that a {\em Borel space} is a set $X$ together with a $\sigma$-algebra of sets $\mathcal B$ (the ``Borel'' sets), and that a function between such spaces is Borel if inverse images of Borel sets are Borel. $X$ is said to be \emph{standard} if either it is countable with all sets Borel, or it is Borel isomorphic with $[0,1]$ with the usual Borel sets.

Mackey was the first to recognize that there is a fundamental obstruction to
generalizing (2.2). As before, we let $\hat{G}$ be the unitary equivalence
classes of (separable) unitary irreducible representations of $G$. In order
to define direct integrals over $\hat{G},$ it is necessary to single out a
natural Borel space structure on $\hat{G}$. Following Mackey, we let $H_{n}$ be the $n$ dimensional Hilbert space $\Bbb{C}^{n}$ for $n\in \Bbb{N}$, and $H_{\infty}=\ell_{2}$, and we let Irr$_{n}\,G$ be the set of all irreducible representations of $G$ on $H_{n}$. We have a natural group
action of the unitary group $\mathcal{U}_{n}$ of $H_{n}$ on Irr$_{n}\,G$ determined by
\[
(U,\pi)\mapsto U^{-1}\pi U
\]
Letting $\hat{G}%
_{n}$ be the orbit space of this action, we define $\hat{G}=\sqcup \hat{G}%
_{n}$ to be the \emph{dual} of $G.$ The difficulty is that although each
space Irr$_{n}\,G$ is a standard Borel space, the \emph{quotient }Borel
structure on $\hat{G}_{\infty }$, and thus that on $\hat{G},$ can be
``non-smooth''. To be more precise, there is a fundamental \emph{dichotomy:}
either $\hat{G}$ is standard, or one cannot distinguish the points of $\hat{G}
$ with countably many Borel functions (``parameters'') $f_{k}:\hat{G}%
\rightarrow \Bbb{R} \,\,(k=1,2,\ldots)$

Mackey apparently believed that if $\hat{G}$ is not smooth, then the
irreducible representations are \emph{unclassifiable. }This phenomenon is
perhaps best illustrated by considering $G=\Bbb{F}_{2}.$ In order to
describe $\widehat{\Bbb{F}_{2}},$ it would be neccesary to classify all
irreducible pairs of unitaries $(U,V)$ to within unitary equivalence. I think that most would regard this as an ``impossible'' task. 

Under favorable circumstances (one says that the group is type I), one need not use Murray and von Neumann's continuous dimensions. In that case all separable representations have canonical direct integral decompositions
\[
\pi \simeq \int_{\hat{G}}^{\oplus}n_{\theta}\theta d\mu(\theta),
\]
where $n_{\theta}\in \Bbb{N}\cup\{0\}\cup\{\infty\}$. 

In the non-type I case, all hell breaks loose. Even if one uses real number multiplicities $n_{\theta}$, the decomposition fails. But the pathology of the situation becomes more apparent when one turns to multiplicity free decompositions, i.e., those of the form 
\[
\pi\simeq\int_{\hat{G}}^{\oplus }\theta
d\mu (\theta )
\]
for ``smooth'' measures $\mu $ on $\hat{G}.$ For example, there are examples of groups $G$ for which there exist \emph{disjoint} measures $\mu$ and $\nu $ on $\hat{G}$ for which the multiplicity-free representations $\int_{\hat{G}}^{\oplus }\theta d\mu
(\theta )$ and $\int_{\hat{G}}^{\oplus }\theta d\nu (\theta )$ are unitarily
equivalent (see \cite{M1}, Ch. 3, \S 5). The bottom line is that \emph{Mautner's decomposition does not reduce the classification of unitary representations to the classification of the irreducible representations}.

From the evidence of this and other examples, Mackey conjectured ``A group has smooth
dual if and only if all of its representations are of type I'' \cite{M1}, Ch. 2, 
\S 2. This was finally proved by James Glimm \cite{G1}, who considered the
(equivalent) task of classifying the space $\hat{A}$ of all unitary
equivalence classes of $*$-representations of a separable C$^{*}$-algebra $A.$ In a truly remarkable paper he showed that the following conditions are equivalent:
\begin{enumerate}
\item[a1)] $\hat{A}_{\infty }=\mathrm{Irr}_{\infty}\,A/\mathcal{U}_{\infty }$ is a standard Borel space

\item[a2)] $A$ has only type I representations (these involve only integer dimensions)

\item[a3)] $A$ is a GCR algebra in the sense of Kaplansky \cite{Kap}.
\end{enumerate}
Roughly speaking, a GCR algebra is a C$^{*}$-algebra constructed by taking extensions of algebras that are built from fields of compact operators.

In a subsequent article \cite{G2}, Glimm showed that there is an analogous dichotomy for locally
compact transformation groups $G\curvearrowright X$. Among other things, he proved
that the following are equivalent:
\begin{enumerate}
\item[b1)] $X/G$ is standard
\item[b2)]$X$ has no non-trivial ergodic measures
\item[b3)]$X/G$ is T$_{0}$
\item[b4)]The orbits $Gx$ are locally closed (or equivalently, locally compact).
\end{enumerate}

It occurred to me that if one could generalize the transformation group result to appropriate Polish
transformation groups such as 
$\mathcal{U}_{\infty }\curvearrowright\mathrm{Irr}_{\infty}\,G$, 
Mackey's conjecture would immediately follow. In particular, if $\mu $ is a non-trivial ergodic measure on $X_{\infty}=$Irr$_{\infty}\,G,$
then $\pi =\int_{X_{\infty}}^{\oplus }\theta d\mu (\theta )$ is a non-type I
factor representation. I accomplished this in \cite{E1} by replacing Glimm's locally compact arguments with first and second category tricks. Subsequently I streamlined
the theory in \cite{E5}. In that form I was able to reprove Mackey's conjecture by using the unitary group of a unital C$^{*}$-algebra acting on the pure state space via inner automorphisms. 

I was aware that Polish transformation groups arise in many classification
problems. It seemed likely that other intractable
classfication problems are in fact ``unsolvable'' in Mackey's sense. In support of this notion, I
considered some ``toy'' examples in \cite{E5}. 

\section{Early results}

It was Mackey's fault that soon after I started working with him, I became a
renegade student. He suggested that I read Kadison's profound operator
algebraic approach to the classification of representations \cite{Kd}. I don't think he could have foreseen that I would be seduced by operator algebra theory. I was captivated
by the notion that operator algebras were not introduced by von Neumann as a vehicle for ``mathematizing'' quantum mechanics. He had already done \emph{that} in \cite{vn1}. Rather he regarded them as a framework for ``quantizing'' mathematics. 

From the very beginning, classification has been a focal problem of operator algebra theory. Von Neumann
had proved that any ``ring of operators'' (a.k.a. von Neumann algebra) $R$
has a direct integral decompostion into factors (von Neumann algebras
with scalar centers): 
\begin{eqnarray}
R\simeq \int_{X}^{\oplus }R(t)d\mu (t)
\end{eqnarray}
(see \cite{vn}). It was thus natural to define a canonical index space $\widehat{\mathcal{F}}$ of all spatial isomorphism classes of von Neumann algebras on a separable Hilbert space,
and to consider corresponding direct integrals over $\widehat{\mathcal{F}}$: 
\[
R\simeq \int_{\widehat{\mathcal{F}}}^{\oplus }Sd\mu (S).
\]
Despite the fact that only finitely many factors had been
discovered at that time, I had hoped one could use indirect methods to show $%
\widehat{\mathcal{F}}$ is non-smooth, and thus, presumably, that there are uncountably many factors. 

As a first step,
I succeeded in proving that the set $\mathcal{F}$ of all factors on a separable Hilbert space $H$ is a
standard Borel space \cite{E2}. Letting $\mathcal{U}$ be the unitary group
on $H,$ the unitary equivalence action $\mathcal{U}\curvearrowright \mathcal{%
F}$ is a Borel transformation group. In this context, the existence of a
non-trivial ergodic Borel measures on $\mathcal{F}$ would imply the existence of \emph{globally
pathological} von Neumann algebras. I also showed that the non-smoothness of the
unitary equivalence classes of factors was equivalent to the non-smoothness of the $*$-isomorphism classes of factors \cite{E3}. However, more difficulties had to be overcome. It was
not apparent to me that $\mathcal{U}\curvearrowright \mathcal{F}$ was a
Polish tranformation group. This was subsequently proved by Marechal (\cite{Mar}, see
also \cite{H1}, \cite{H2}). The next step would have been to show that there exist
orbits $\mathcal{O}_{1}$ and $\mathcal{O}_{2}$ which are mutually entangled, i.e.,
the closure of each orbit contains the other, since then the quotient space would not be T$_{0}$.

However the program was
soon overtaken by more constructive arguments. Using the multiplicity of
hyperfinite type III von Neumann algebras discovered by Powers \cite{P}, Nielsen
showed that the space $\hat{\mathcal{F}}$ is non-smooth, and in fact there exist non-trivial ``global factors'' \cite{N}. \emph{Thus the idea that central decompositions (3.1) reduce classification theory for von Neumann algebras to that for factors is simply not the case.} Despite the fact that my original goal was frustrated, the result that $\mathcal{F}$ is standard greatly simplifies von Neumann's direct integral theory (see \cite{N1}). 

A few years later, an application for a Polish transformation group
dichotomy unexpectedly arose in the study of homogenous spaces. If $X$ is
a compact metric space and $G$ is its group of homeomorphisms, then $G\curvearrowright X$ is a 
Polish transformation group \cite{Hur}. In particular, if $X$ is homogeneous, i.e., the action is transitive, then $\hat{X}=X/G$ consists
of just one point, and thus it is smooth. Surprisingly, the new information provided by some of the conditions equivalent to smoothness furnished the solution to a problem of some importance in the subject (see \cite{R} and \cite{Macias} for references). In the general
(non-homogeneous) case, one may regard $X/G$ as an index space for a ``homogeneous classification of the points'' in $X.$

\section{Algebraic problems}

As a closet algebraist, I have often been attracted to
algebraic classification problems. At an early point in my graduate career, I was enchanted by Kaplansky's elementary monograph on the classification problems for countable abelian groups \cite{Kp}. Whereas the finitely generated abelian groups are usually classified in undergraduate modern algebra
courses, the situation for non-finitely generated abelian groups is quite
mysterious. This naturally seemed to me to be a promising area to test Mackey's
classification philosophy. My investigations in this direction were stymied
by a series of difficulties.

The countable torsion abelian groups are completely classified by the Ulm
invariants. This did not appear to support Mackey's philosophy,
since one must use uncountably many invariants. We will briefly revisit this case below.

Turning to the countable \emph{torsion free} abelian groups, the rank one
groups, i.e., those which can be embedded in $\Bbb{Q},$ were classified by Baer \cite{B}. Sketching his results, we define a \emph{generalized natural number} $m$ to be a mapping $m:\Bbb{P}\to{0,1,2,\ldots,\infty}$ where $\Bbb{P}=\{2,3,\ldots\}$ is the set of primes. It is convenient to use the notation
\[
m= 2^{a_{1}}3^{a_{2}}\ldots
\]
where $a_{i}=m(p_{i})$. We let $\mathfrak{N}$ be the set of generalized integers, and we identify $\Bbb{ N}$ with the $m\in \mathfrak{N}$ for which $\sum a_{i}<\infty.$ 

Each $m\in \mathfrak{N}$ determines a subgroup  $G(m)$ of $\Bbb Q$ via 
\[ G(m)=\{r/s:r\in \Bbb{Z}, s\in \Bbb {N},s|m\},
\]
where we use the obvious notion of divisor. All rank one groups arise in this fashion. Furthermore if $n= 2^{b_{1}}3^{b_{2}}\ldots$, then $G(m)\simeq G(n)$ if and only if $a_{k}=b_{k}$ for all but finitely many $k$, and $a_{k}=\infty$ if and only if $b_{k}=\infty$. We have, for example, that $G(n)\simeq \Bbb{Z}$ if and only if $n\in \Bbb{N}$, whereas $G(2^{\infty})=\Bbb{Z}[1/2]$. Again this is regarded as a completely
satisfactory classification, despite the fact that it is not smooth (see below).

The classification of the rank $n\geq 2$ groups (those which can be embedded in $%
\Bbb{Q}^{n}$ but not in $\Bbb{Q}^{n-1})$ has long been recognized to be
intractable (see \cite{F}). One has examples of non-isomorphic finite rank groups 
$G_{1},G_{2},$ $H_{1},$ $H_{2}$ which are not themselves decomposable into
direct sums, for which 
\[
G_{1}\oplus G_{2}\simeq H_{1}\oplus H_{2}. 
\]

More generally, one can consider the classification of arbitrary countable groups by considering the set $\mathcal{G}$ of all
mappings 
\[
m:\Bbb{N}\times \Bbb{N}\rightarrow \Bbb{N}:(s,t)\rightarrow s\cdot t 
\]
which satisfy the usual group actions. The group $S_{\infty }$ of all
permutations of $\Bbb{N}$ acts on $\mathcal{G}$ via $\pi (m)(s,t)=\pi^{-1}m(\pi(s),\pi(t)).$ If we let $S_{\infty }$ and $\mathcal{G}
$ have the relative topologies in $\Bbb{N}^{\Bbb N},$ and $\Bbb{N}^{\Bbb{N}\times 
\Bbb{N}},$ respectively, we obtain a Polish transformation group $S_{\infty
}\curvearrowright \mathcal{G},$ and we may regard the quotient space $%
\widehat{\mathcal{G}}$ as the ``classfication space'' for countable discrete
groups. Unfortunately, I could see no route to proving results about this space. Furthermore,
the use of transformation groups seemed awkward. In principle one should be
able to find dichotomies for suitable equivalence relations on Borel spaces.

Owing in part to my limited knowledge of descriptive set theory (my only
source being Kuratowski's book \cite{Ku}), I could not push the program any further. My subsequent move from the University of Pennsylvania to UCLA in 1979 was quite
fortuitous, since Los Angeles has long been a world class center for the
study of mathematical logic, owing to the efforts of Yiannis Moschovakis, Tony Martin, and Aleko Kechris. It has been a delight to witness the
establishment of Borel classification as a vital area of mathematics due the work of such individuals as Adams, Harrington, Hjorth, Kechris, Louveau, and Thomas and a 
host of others.

\section{Enter the logicians}

Although I had approached a few logicians at an earlier point, I could not
have anticipated the sudden explosion of results that occured in the 1990's.
The descriptive set theorists formulated the Borel classification problems
in terms of ``Borel cardinality'' for equivalence relations. Given a
standard Borel space $X$, an equivalence relation $E$ on $X$ is said to be
Borel if $E$ is a Borel subset of $X\times X.$ We let $X/E$ denote the
equivalence classes with the quotient Borel structure, i.e. $B\subseteq X/E$
is said to be Borel if and only if its inverse image in $X$ is Borel. We
simply write $X/E=X$ if $E$ is the trivial equivalence relation of equality.
Given two equivalence relations $E$ and $F$ on standard Borel spaces $X$ and 
$Y,$ we write $E\leq _{B}F$ (and say that $E$ has smaller Borel cardinality
than $F)$ if there is a Borel map $f:X\rightarrow Y$ such that $xEy$ if and
only if $f(x)Ff(y).$ In the language of \cite{Hj2} (\emph{from which this discussion
is purloined} -- I am out of my depth here), there is a ``Borel embedding'' of $X/E$ into $Y/F$ in the sense
that the injection has a Borel lifting. We let $E\sim_{B}F$ if $E\leq _{B}F$ and $F\leq _{B}E$, and $E<_{B}F$ if $E\leq _{B}F$ but the reverse is not true. Of particular importance is the Vitali
relation $E_{0}$ on $2^{\Bbb{N}}$ of sequences $(a_{n}),$ $a_{n}\in \left\{
0,1\right\} $ where $aE_{0}b$ if and only if $a_{n}=b_{n}$ for all but
finitely many $n.$ One has that $E_{0}\sim E_{V},$ where $E_{V}$ is
the Vitali equivalence on $\Bbb{R}$ determined by $xE_{V}y$ if and only if $%
y-x\in \Bbb{Q}.$

An early ``level 0'' dichotomy for an arbitrary Borel (or even coanalytic) equivalence relation $E
$ was discovered by Silver \cite{S} (see also \cite{Be}):
\[
E\leq _{B}\Bbb{N}\text{ or }\Bbb{R}\leq _{B}E.\text{ }
\]
In particular, if $X/E$ is not countable, then $X$ has a perfect subset of
pairwise inequivalent elements.

In a truly striking application of computability theory techniques, that should probably
be in the tool-kit of any serious functional analyst who uses measure theory, Harrington, Kechris, and
Louveau proved what seems to be the ``best possible theorem'' at ``level 1'':
\[
E\leq _{B}\Bbb{R}\text{ or }E_{0}\leq _{B}E.
\]
This is a far-reaching generalization of the methods of Glimm and myself
since at the heart of our arguments was a construction of an embedding of $%
E_{0}$ into the relevant quotient spaces. A total of seven dichotomies is
surveyed in \cite{Hj2}, which delineate a remarkable non-linear architecture among
the ''Borel cardinals''. 

But it is in the realm of applications that some of the most breathtaking
results have been proved. Returning to the question of why the finite
torsion free groups are intractable, let $S=S(\Bbb{Q}^{n})$ be the Polish
space of additive subgroups of $\Bbb{Q}^{n},$ and $\simeq _{n}$ denote the
algebraic isomorphism equivalence on $S$. The latter is  the orbit space of the natural countable group action $GL_{n}(\Bbb{Q})\curvearrowright S.$ Following up work of Adams, and Kechris
and Hjorth, who utilized Zimmer's superrigidity theory (!), Simon Thomas proved
that one has the hierarchy
\begin{eqnarray}
E_{0}\,\sim _{B}\,\simeq _{1}\,<_{B}\,\simeq _{2}\,<_{B}\ldots 
\end{eqnarray}
(see \cite{Th}). The Borel equivalence in this chain just reflects the solution of
the classification problem for rank one groups. (5.1) provides a genuine
descriptive set-theoretic explanation of the difficulty algebraists have had
with classifying these groups.  

In a remarkable twist, the general theory is in fact related to Sorin Popa's
theory of superrigidity theory for factors \cite{Pop}. But we are in no position to go
into that!

\section{Some conjectures}

Classification problems pervade all branches of mathematics, and each time
an intractable situation arises, the possibility of Borel cardinal obstructions to classification should be considered. Thus it now seems unlikely that anyone will try to classify the subgroups of $\Bbb Q^{2}$. I will survey some intriguing
questions that remain unanswered at this point. 

It might appear that Mackey was too swift in proscribing the non-smooth
equivalence relations as being unclassifiable. The rank one torsion groups
have a non-smooth isomorphism relation (see above), but these are regarded as ``classified'' (see \cite{Hj} for a discussion of the Ulm invariants). In particular
we might well define a Borel equivalence relation $E$ to be
``unclassifiable'' if $E_{0}<_{B}E.$ On the other hand, perhaps Mackey was
correct. To my knowledge, there is no non-type I group $G$ for which $\hat{G}$ has been
classified in any reasonable sense. This leads us to refine Mackey's conjecture to state that if $G$
is a non-type I group, then $\hat{G}>_{B}E_{0}$, i.e., ``non-type I if and only if very non-smooth dual''. We would also expect that there will be a hierarchy for the Borel cardinals of group duals. For those seeking
non-type I groups with which to experiment, we recall Thoma's theorem that
the only countable type I discrete groups are those that are finite
extensions of abelian groups \cite{T}. 

After reading the above, Aleko Kechris kindly informed me that indeed Greg Hjorth has shown that non-type I discrete groups have a \emph{very} non-smooth dual. An account of that result may be found in Appendix H of \cite{KN}. A great many other applications of classification are also considered.

But there are many intractable problems for which the Borel
cardinals seem less useful. My favorite example is the task of
classifying nilpotent Lie algebras of a given dimension. There are exactly
4, 231, 9,022, 1,028, and 10,285 nilpotent Lie algebras of dimensions 4, 5, 6, 7, and 8, respectively, over a field of characteristic 0 (see \cite{Ts}, I am quoting from the AMS review by Oleg Gutik). On the other hand, I. Magnin \cite{Mag} wrote that ``La
classification des alg\`{e}bres de Lie nilpotentes de dimension finie est un
probl\`{e}me non r\'{e}solu \`{a} ce jour ...'' An evidence of pathology at
higher dimensions and the real field can be found in a paper by Chong-yun
Chao \cite{C}, who by selecting real structure constants that are linearly
independent over $\Bbb{Q}$, showed that there are uncountably many
nonisomorphic nilpotent Lie algebras of dimension 10 and length 2. Could this be another instance of non-smoothness?

If we let $\mathcal{L}_{n}$ be the set of
\emph{all} Lie algebras of dimension $n$, the isomorphism
relation is determined by the obvious locally compact transformation group  $GL(n,\Bbb{R})\curvearrowright \mathcal{L}_{n}$. But in fact this is an algebraic
transformation group, and from a theorem of Chevalley (see \cite{D}, p. 183, bottom, and \cite{Hum}, page 60) the orbits are necessarily locally closed. Thus
the quotient space is smooth and we must declare that \emph{all} the Lie algebras of dimension $n$ are classifiable.

An asymptotic approach might yield a more interesting result. Consider the problem of classifying all of the finite dimensional Lie algebras ``at once''. We define $\mathcal{L}=\cup \mathcal{L}_{n}$. In this context there are natural mappings 
\[
\mathcal{L}_{n}\hookrightarrow \mathcal{L}_{n+m}:L\to L\oplus 0.
\]
If we broaden the notion of isomorphism to include such links between the spaces
$\mathcal{L}_{n}$, it is possible that the resulting quotient space would not simply be a disjoint union of smooth spaces, and might even be non-smooth.

A much more speculative question is whether limiting arguments could be used to study the classification of finite groups. Turning this around, one might be able to show that finitely presented groups have a non-smooth classification, which somehow ``infiltrates'' the finite group problem. 

In such a context one might consider rings of sets rather than $\sigma$-rings. 
Perhaps the situation is analogous to the the theory of finitely additive measures on normal spaces, where one uses the ring generated by the locally closed sets. The finitely additive measures correspond to the countably additive Radon measures on the $\beta$-compactification of the normal space. Presumably one would want to use a smaller, metrizable compactification. A. Kechris has pointed out to me that computably enumerable equivalence relations on a countable set have been considered in \cite{GG}.

Other problems that I have encountered that might also be approached in this spirit are the order-isomorphic classification of finite rank dimension groups \cite{E4} and the completely isometric classification of finite dimensional operator spaces \cite{E6}.

The asymptotic notion might also be related to questions of computability. Going very far out on a limb, could the notorious $P$, $NP$ problem have something to do with a finitistic Mackey dichotomy?

\section{A final classfication problem}

As a shameless proselyte of mathematical quantization (see, e.g. \cite{E7}), I have long puzzled over how I might have converted George to the ``cause''. Over the years, I have had many imaginary conversations with him on the subject. Murray and von Neumann's articles can be regarded as the quantization of integration theory. The depth of the theory is reflected by the fact that it took nearly fifty years before it was completed, owing to the pioneering work of Tomita and Takesaki \cite{T} (see \cite{KC}). In the last few decades, we have witnessed the development of quantized geometry, largely due to Alain Connes (see \cite{Co}). To my astonishment, my own favorite quantization (quantized Banach space theory or ``operator spaces'', see \cite{E6}, \cite{Pis}) is currently being used in a highly non-trivial fashion by mathematicians, physicists and engineers studying entropy for quantum channels (see, e.g., \cite{Ju}, and \cite{Choi}).

I knew that George never fully subscribed to the quantization of mathematics begun by von Neumann. At a dinner arranged by Irving Kaplansky and his wife during the operator algebra year (2000-2001) at MSRI, I had the pleasure of introducing my wife Rita to Alice and George. During that occasion, I asked him what he thought about the theory of quantum groups. I suppose that I was still trying to explain to him why I had abandoned group representation theory some forty years before. He answered that he was not happy with the terminology. I pointed out that mathematicians often use terms such as ``dynamical systems'' in non-physical situatons. He demurred, saying that we should discuss it at some other time. We never had that oppurtunity. 
	
The delay between the appearance of the Murray and von Neumann's papers and the first recognition of their significance by Mackey's young colleagues Segal, and then Kadison, Kaplansky and Singer, was due to many factors. To begin with, the notions of \emph{classical} functional analysis  were still quite novel, and few were willing to tackle the quantum theoretic analogues introduced by von Neumann. Furthermore, there was a lingering hope by many mathematicians that quantum theory and its fiendishly nonclassical ``paradoxes'' might be a transitory phenomenon that could be ignored. Of course, the Second World War was detrimental to abstract research. But it must also be said that von Neumann had a very difficult writing style, in which he ``took no prisoners''.  He apparently felt that the material spoke for itself, and that it did not need motivation. 
	
In the manner of many brilliant investigators, Mackey preferred coming up with his own ideas to reading the work of others. It is my understanding that he reconstructed the Murray and von Neumann theory of types before he realized that it had already been developed in the Rings of Operator papers. Nevertheless, in contrast with most of the functional analytical community in the United States during the fifties and sixties, Mackey fully understood the importance to mathematicians of understanding the foundations of quantum physics. His monograph on the subject \cite{M3} is a landmark in the mathematical literature. Thus my conclusion is that George was a transitional figure. 
  
But I think that he would have disagreed. In fact his interest in so many branches of mathematics made him \emph{unclassifiable}. 
He was above all else a generalist who sought to break down the walls that are imposed by our specialized backgrounds. I am indebted to him for exemplifying that philosophy to me.

\end{document}